\documentclass[oneside]{article}
\usepackage{sibgrapi}
\usepackage{times}
\usepackage{amssymb}
\usepackage{amsmath}
\usepackage{subfigure}
\usepackage{graphicx}
\usepackage{subfigure}
\begin{document}

\newtheorem{theorem}{Theorem}
\newtheorem{corollary}[theorem]{Corollary}
\newtheorem{lemma}[theorem]{Lemma}
\newtheorem{proposition}[theorem]{Proposition}
\newtheorem{definition}[theorem]{Definition}
\newtheorem{example}{Example}
\newtheorem{conjecture}[theorem]{Conjecture}
\newenvironment{remark}[1][Remark.]
{\begin{trivlist} \item[\hskip \labelsep {\bf #1}]}{\end{trivlist}}
\newenvironment{proof}[1][Proof.]
{\begin{trivlist} \item[\hskip \labelsep {\it \bfseries
#1}]}{\hfill $\blacksquare$ \end{trivlist}}

\newcommand{\Z}[1][]{\ensuremath{{\mathbb{Z}^{#1}} }}
\newcommand{\C}[1][]{\ensuremath{{\mathbb{C}^{#1}} }}
\newcommand{\R}[1][]{\ensuremath{{\mathbb{R}^{#1}} }}
\renewcommand{\H}[1][]{\ensuremath{{\mathbb{H}^{#1}} }}
\renewcommand{\S}[1][]{\ensuremath{{\mathbb{S}^{#1}} }}
\newcommand{\<}{\langle}
\renewcommand{\>}{\rangle}
\newcommand{\eps}{\epsilon}

\renewcommand{\v}{\mathbf{v}}
\renewcommand{\P}{\mathbf{P}}
\newcommand{\p}{\mathbf{p}}
\newcommand{\q}{\mathbf{q}}
\newcommand{\w}{\mathbf{w}}
\newcommand{\x}{\mathbf{x}}
\newcommand{\y}{\mathbf{y}}
\newcommand{\zb}{\overline{z}}
\renewcommand{\t}{\mathbf{t}}
\newcommand{\n}{\mathbf{n}}
\renewcommand{\b}{\mathbf{b}}
\newcommand{\T}{\mathbf{T}}
\renewcommand{\L}{\mathbf{L}}
\newcommand{\V}{\mathbf{V}}
\newcommand{\X}{\mathbf{X}}
\newcommand{\Y}{\mathbf{Y}}
\newcommand{\N}{\mathbf{N}}
\newcommand{\dl}{\Delta l}
\renewcommand{\k}{\kappa}
\newcommand{\norm}[1]{\Vert #1 \Vert}
\newcommand{\fsep}{\hspace*{\fill}}



\title{ Discrete Affine Minimal Surfaces with Indefinite Metric}

\author{Marcos Craizer\authortag{1}
  Henri Anciaux\authortag{2},
  and
  Thomas Lewiner\authortag{1}
} 

\address{
\authortag{1} Department of Mathematics, Pontif\'{i}cia Universidade Cat\'{o}lica, \\
 Rio de Janeiro, Brazil.\\
 \authortag{2}Department of Mathematics and Computing, Institute of Technology, Tralee, \\
 Co. Kerry, Ireland. \\
{\small \tt craizer@mat.puc-rio.br},\ \ {\small \tt
henri.anciaux@staff.ittralee.ie},\ \
 {\small \tt tomlew@mat.puc-rio.br} \\
}

\abstract {Inspired by the Weierstrass representation of smooth
affine minimal surfaces with indefinite metric, we propose a
constructive process  producing a large class of discrete surfaces
that we call \em discrete affine minimal surfaces. \em We show
that they are  critical points of an affine area functional
defined on the space of quadrangular discrete surfaces. The
construction makes use of asymptotic coordinates and allows
defining the discrete analogs of some differential geometric
objects, such as the normal and co normal vector fields, the cubic
form and the compatibility equations.
\\[2mm] \noindent
 }

\maketitle

\section{Introduction}

In  affine differential geometry, the notion of minimal surfaces,
i.e. the critical points of the affine area functional, arises
naturally and has received a broad attention in the last decades.
In particular, it has been proved in
\cite{Calabi82},\cite{Calabi88} that convex affine minimal
surfaces actually maximize the affine area, thus justifying the
sometimes used terminology {\it maximal surfaces}. On the other
hand, \cite{Vrancken89} showed that this is not true for
non-convex surfaces. In the convex or non-convex case,
Weierstrass-type representations have been derived, allowing the
explicit construction of local parameterizations of affine minimal
surfaces from the co-normal vector field. This representation
makes use of isothermal coordinates in the definite case and
asymptotic coordinates in the indefinite case.

More recently, the expansion of computer graphics and applications
in mathematical physics have given a great impulse to the issue of
giving discrete equivalents of differential geometric objects
(\cite{Bobenko06},\cite{Bobenko08}). In the particular case of
affine geometry some work has been done toward a theory of
discrete affine surfaces. In \cite{Bobenko199} a consistent
definition of discrete affine spheres is proposed, both for
definite and indefinite metrics and in \cite{Matsuura03} a similar
construction is done in the context of improper affine spheres.

 In this work we introduce a discrete
analog of the smooth Weierstrass representation in the indefinite
case, giving rise to explicit parameterizations of quadrangular
surfaces in discrete asymptotic coordinates that we call {\it
discrete affine minimal surfaces}. Over these discrete affine
minimal surfaces, we can define the {\it discrete affine metric},
the {\it discrete affine normal vector field} and a discrete
analog of the smooth cubic form, that we shall call {\it discrete
affine cubic form}. We show that, as occurs in the smooth case,
the discrete affine metric and the discrete affine cubic form must
satisfy compatibility equations. Moreover, these compatibility
equations are a necessary and sufficient condition for the
existence of an affine minimal surface, given its metric and cubic
form.

We also introduce a natural affine area functional in the set of
quadrangular indefinite discrete surfaces and show that the
minimal surfaces that we have constructed are critical points of
this functional, thus justifying the choice of our terminology.

In view of the above results, it is natural to ask wether it is
possible drop the minimality condition in this construction. This
issue is related to the problem of finding a convenient definition
of discrete affine mean curvature vector. In another direction, it
is tempting to look for an analogous construction in the definite
case. We plan to address these questions in a forthcoming work.

The paper is organized as follows: in Section 2 we state some
classical notations and facts about asymptotic parameterizations
of indefinite affine smooth surfaces in $\mathbb R^3.$ In Section
3, inspired by the continuous case, we implement the construction
process of discrete affine minimal surfaces. Section 4 is devoted
to the description of the variational property of these surfaces
(Theorem \ref{variational}). In last section, we introduce the
discrete affine cubic form, derive the compatibility equations and
prove the corresponding theorem of existence and uniqueness
(Theorem \ref{Existence}).


\section{Preliminaries}

\smallskip\noindent{\bf Notation.} Along the paper, letters in subscripts denote partial
derivatives with respect to the corresponding variable, and
$V_1\cdot V_2$, $[V_1,V_2,V_3]$ and $V_1\times V_2$ denote
respectively the inner product, the determinant and the
cross-product of vectors $V_1,V_2,V_3\in\R^3$.

\medskip

Consider a parameterized smooth surface $q:U\subset \R^2\to \R^3$,
where $U$ is an open subset of the plane and denote by
\begin{eqnarray*}
L(u,v)&=&[q_u,q_v,q_{uu}]\\
M(u,v)&=&[q_u,q_v,q_{uv}]\\
N(u,v)&=&[q_u,q_v,q_{vv}]\\
\end{eqnarray*}
The surface is non-degenerate if $LN-M^2\neq 0$, and, in this
case, the Berwald-Blaschke metric is defined by
$$
ds^2=\frac{1}{|LN-M^2|^{1/4}}\left({Ldu^2+2Mdudv+Ndv^2}\right)
$$
If $LN-M^2>0$, the metric is {\it definite} while if $LN-M^2<0$,
the metric is {\it indefinite}. In this paper, we shall restrict
ourselves to surfaces with indefinite metric.

We say that the coordinates $(u,v)$ are {\it asymptotic} if
$L=N=0$. In this case, the metric takes the form $ds^2=2Fdudv$,
where $F^2=M$. Also, we can write
\begin{eqnarray}
q_{uu}&=&\frac{1}{F}(F_uq_u+Aq_v)\label{struct1Smootha}\\
q_{vv}&=&\frac{1}{F}(Bq_u+F_vq_v) \label{struct1Smoothb},
\end{eqnarray}
where $A=A(u,v)$ and $B=B(u,v)$ are the coefficients of the affine
cubic form $Adu^3+Bdv^3$ (see \cite{Nomizu94}).

The vector field $\xi(u,v)=\displaystyle{\frac{q_{uv}}{F}}$ is
called the {\it affine normal} vector field. We have
\begin{eqnarray}
\xi_{u}&=&-Hq_u+\frac{A_v}{F^2}q_v\label{struct2Smootha}\\
\xi_{v}&=&\frac{B_u}{F^2}q_u-Hq_v\label{struct2Smoothb},
\end{eqnarray}
where $H$ is the {\it affine mean curvature}. Equations
\eqref{struct1Smootha}, \eqref{struct1Smoothb},
\eqref{struct2Smootha}  and \eqref{struct2Smoothb} are the
structural equations of the surface. For a given surface, the
quadratic form $Fdudv$, the cubic form $Adu^3+Bdv^3$ and the
affine mean curvature $H$ should satisfy the following
compatibility equations:

\begin{eqnarray}
H_{u}&=&\frac{AB_u}{F^3}-\frac{1}{F}(\frac{A_v}{F})_v,\label{CompatibilitySmootha}\\
H_{v}&=&\frac{BA_v}{F^3}-\frac{1}{F}(\frac{B_u}{F})_u.
\label{CompatibilitySmoothb}
\end{eqnarray}
Conversely, given $F,A,B$ and $H$ satisfying equations
\eqref{CompatibilitySmootha} and \eqref{CompatibilitySmoothb},
there exists a parameterization $q(u,v)$ of a surface with
quadratic form $2Fdudv$, cubic form $Adu^3+Bdv^3$ and affine mean
curvature $H$. For details of the above equations, see
\cite{Buchin83}.

The vector field $\nu(u,v)=\frac{q_u\times q_v}{F}$ is called the
{\it co-normal} vector field. It satisfies Lelieuvre's equations
\begin{eqnarray}
q_u&=&\nu\times\nu_u
\label{LelieuvreSmootha}\\q_v&=&-\nu\times\nu_v.\label{LelieuvreSmoothb}
\end{eqnarray}
 It also satisfies the equation
$\Delta\nu=-2H\nu$ , where $\Delta$ denotes the Laplacian with
respect to the Berwald-Blaschke metric (e.g., see
\cite{Nomizu94}). It turns out that in asymptotic coordinates,
$\Delta\nu= \nu_{uv}$.

A surface is said to be {\it affine minimal } if its affine mean
curvature $H$ vanishes or equivalently if its co-normal vector
field satisfies the equation $\nu_{uv}=0$. The interest of the
co-normal definition lies in the fact that the resolution of this
last equation is straightforward: $\nu_{uv}=0$ if and only if
$\nu(u,v)$ takes the form $\nu(u,v)=\nu^1(u)+\nu^2(v)$, where
$\nu^1$ and $\nu^2$ are two real functions of one variable.
Starting from the co-normal vector field and using Lelieuvre's
equations \eqref{LelieuvreSmootha} and \eqref{LelieuvreSmoothb},
one gets an immersion $q$ which turns to be a parameterization in
asymptotic coordinates of an affine minimal surfaces. This is a
simple way to construct examples of smooth affine minimal surfaces
(e.g., see \cite{Simon93}).






\section{Definitions, properties and
examples}\label{SectionDefinitions}

In this section, inspired by the properties of affine minimal
surfaces and asymptotic coordinates discussed above, we describe a
construction process of a class of discrete surfaces with
properties analogous to the smooth case. We start
 with a vector field of the form
$\nu(u,v)=\nu^1(u)+\nu^2(v),$ where $\nu^1$ and $\nu^2$ are two
real functions of one discrete variable. In particular $\nu$ is
the restriction to a subset of $\Z^2$ of a smooth co-normal vector
field of some smooth minimal surface. To obtain the affine
immersion, we make a discrete integration of the discrete analogs
of Lelieuvre's equations \eqref{LelieuvreSmootha} and
\eqref{LelieuvreSmoothb}.



\smallskip\noindent{\bf Notation.} For a discrete real or vector function $f:D\subset\Z^2$,
we denote the discrete partial derivatives with respect to $u$ or
$v$ by
\begin{eqnarray*}
f_1(u+\tfrac{1}{2},v)&=&f(u+1,v)-f(u,v)\\
f_2(u,v+\tfrac{1}{2})&=&f(u,v+1)-f(u,v).
\end{eqnarray*}
The second order partial derivatives are defined by
\begin{eqnarray*}
f_{11}(u,v)&=&f(u+1,v)-2f(u,v)+f(u-1,v)\\
f_{22}(u,v)&=&f(u,v+1)-2f(u,v)+f(u,v-1)\\
f_{12}(u,v)&=&f(u+1,v+1)+f(u,v)-f(u+1,v)-f(u,v+1).
\end{eqnarray*}

\subsection{Starting with co-normals}

Consider a map $\nu:D\subset\Z^2\to\R^3$, called the {\it discrete
co-normal map}, satisfying
\begin{equation}\label{HarmonicConormals}
\nu_{12}(u+\tfrac{1}{2},v+\tfrac{1}{2})=0,\ (u,v)\in D.
\end{equation}
We shall also assume that
\begin{equation*}
F(u+\tfrac{1}{2},v+\tfrac{1}{2})=\nu(u,v)\cdot(\nu(u,v+1)\times\nu(u+1,v))>0,\
(u,v)\in D.
\end{equation*}
Discrete co-normal maps can be obtained from smooth maps
$\nu:U\subset\R^2\to\R^3$ satisfying $\nu_{uv}=0$ by restricting
the domain to a subset $D\subset\Z^2$.


\subsection{The affine immersion}

We define the affine immersion by the discrete analog of
Lelieuvre¥s formulas:

\begin{eqnarray}
q_1(u+\tfrac{1}{2},v)&=&\nu(u,v)\times\nu(u+1,v)\label{LeliuvreDiscretasa}\\
q_2(u,v+\tfrac{1}{2})&=&-\nu(u,v)\times\nu(u,v+1).\label{LeliuvreDiscretasb}
\end{eqnarray}

\begin{theorem} \label{DefMinimas} There exists an immersion $q(u,v)$ such that $q_1(u+\tfrac{1}{2},v)$ and
$q_2(u,v+\tfrac{1}{2})$ are as above. Moreover, it satisfies the
following properties:

\begin{enumerate}

\item  The co-normal at $(u,v)$ can be obtained by any of the
following formulas:
\begin{eqnarray*}
\nu(u,v)&=&\frac{1}{F(u+\tfrac{1}{2},v+\tfrac{1}{2})}(q_1(u+\tfrac{1}{2},v)\times
q_2(u,v+\tfrac{1}{2}))\\
\nu(u,v)&=&\frac{1}{F(u-\tfrac{1}{2},v+\tfrac{1}{2})}(q_1(u-\tfrac{1}{2},v)\times
q_2(u,v+\tfrac{1}{2}))\\
\nu(u,v)&=&\frac{1}{F(u-\tfrac{1}{2},v-\tfrac{1}{2})}(q_1(u-\tfrac{1}{2},v)\times
q_2(u,v-\tfrac{1}{2}))\\
\nu(u,v)&=&\frac{1}{F(u+\tfrac{1}{2},v-\tfrac{1}{2})}(q_1(u+\tfrac{1}{2},v)\times
q_2(u,v-\tfrac{1}{2})). \end{eqnarray*}

\item  The parameterization is
asymptotic:

\begin{eqnarray*}
\left[q_1(u\pm\tfrac{1}{2},v),q_2(u,v\pm\tfrac{1}{2}),q_{11}(u,v)\right]&=&0\\
\left[q_1(u\pm\tfrac{1}{2},v),q_2(u,v\pm\tfrac{1}{2}),q_{22}(u,v)\right]&=&0
\end{eqnarray*}
and
\begin{eqnarray*}
\left[q_1(u+\tfrac{1}{2},v),q_2(u,v+\tfrac{1}{2}),q_{12}(u+\tfrac{1}{2},v+\tfrac{1}{2})\right]&=&
F^2(u+\tfrac{1}{2},v+\tfrac{1}{2})\\
\left[q_1(u+\tfrac{1}{2},v),q_2(u+1,v+\tfrac{1}{2}),q_{12}(u+\tfrac{1}{2},v+\tfrac{1}{2})\right]&=&
F^2(u+\tfrac{1}{2},v+\tfrac{1}{2})\\
\left[q_1(u+\tfrac{1}{2},v+1),q_2(u,v+\tfrac{1}{2}),q_{12}(u+\tfrac{1}{2},v+\tfrac{1}{2})\right]&=&
F^2(u+\tfrac{1}{2},v+\tfrac{1}{2})\\
\left[q_1(u+\tfrac{1}{2},v+1),q_2(u+1,v+\tfrac{1}{2}),q_{12}(u+\tfrac{1}{2},v+\tfrac{1}{2})\right]&=&
F^2(u+\tfrac{1}{2},v+\tfrac{1}{2}).
 \end{eqnarray*}

\end{enumerate}
\end{theorem}

\begin{proof}
For the existence of $q$, we must show that the finite difference
equations \eqref{LeliuvreDiscretasa} and
\eqref{LeliuvreDiscretasb} are integrable, i.e.,
$q_{12}-q_{21}=0$. By definition,
\begin{eqnarray*}
q_{12}(u+\tfrac{1}{2},v+\tfrac{1}{2}
)&=&\nu(u,v+1)\times\nu(u+1,v+1)-\nu(u,v)\times\nu(u+1,v)\\
q_{21}(u+\tfrac{1}{2},v+\tfrac{1}{2}
)&=&-\nu(u+1,v)\times\nu(u+1,v+1)+\nu(u,v)\times\nu(u,v+1).
\end{eqnarray*}
Hence
\begin{equation*}
q_{12}-q_{21}=(\nu(u+1,v)+\nu(u,v+1))\times(\nu(u+1,v+1)+\nu(u,v)),
\end{equation*}
which vanishes from property \eqref{HarmonicConormals}.

We now prove only one of the equations of item 1, since the proofs
of the others are similar:
\begin{eqnarray*}
q_1(u+\tfrac{1}{2},v)\times
q_2(u,v+\tfrac{1}{2})&=&-(\nu(u,v)\times\nu(u+1,v))\times(\nu(u,v)\times\nu(u,v+1))\\
&=&-\left[\nu(u,v),\nu(u+1,v)),\nu(u,v+1)\right]\nu(u,v)\\
&=&F(u+\tfrac{1}{2},v+\tfrac{1}{2})\nu(u,v).
\end{eqnarray*}

For the proof of item
 2, we prove one formula of the
 first group and one formula of the second group, the others being
 similar:
\begin{eqnarray*}
L(u+\tfrac{1}{2},v+\tfrac{1}{2})&=&
F(u+\tfrac{1}{2},v+\tfrac{1}{2})\nu(u,v)\cdot(-q_1(u-\tfrac{1}{2},v))\\
&=&F(u+\tfrac{1}{2},v+\tfrac{1}{2})\nu(u,v)\cdot
(\nu(u,v)\times\nu(u-1,v))=0.
\end{eqnarray*}
And
\begin{eqnarray*}
M(u+\tfrac{1}{2},v+\tfrac{1}{2}) &=&
F(u+\tfrac{1}{2},v+\tfrac{1}{2})\nu(u,v)\cdot(q_2(u+1,v+\tfrac{1}{2}))\\
&=& F(u+\tfrac{1}{2},v+\tfrac{1}{2})\nu(u,v)\cdot(\nu(u+1,v+1)\times\nu(u+1,v))\\
&=& F^2(u+\tfrac{1}{2},v+\tfrac{1}{2}),
\end{eqnarray*}
thus completing the proof of the proposition. \end{proof}

The affine immersion $q:D\subset\Z^2\to\R^3$ defined by formulas
\eqref{LeliuvreDiscretasa} and \eqref{LeliuvreDiscretasb} is
called a {\it discrete affine minimal map } and its image a {\it
discrete affine minimal surface }. Along this paper, when there is
no risk of confusion, we shall refer to a discrete affine minimal
map simply as a minimal surface.

\bigskip

\begin{figure}[htb]
\centering  {
\includegraphics[width=0.5
\linewidth,clip =false]{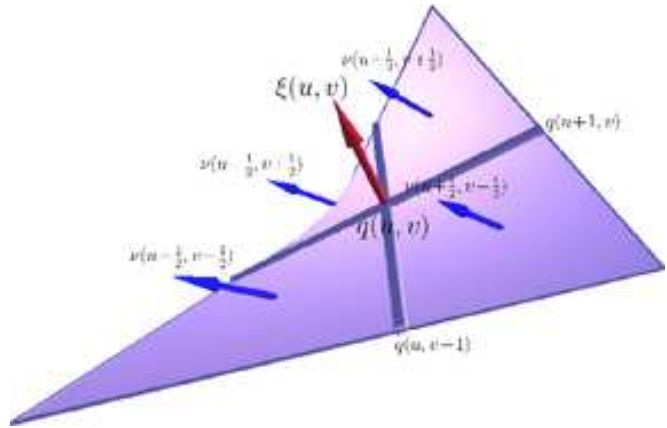}} \caption{The planar cross, the
co-normal vector at the vertex and the normal vectors at the
faces.} \label{figure1}
\end{figure}

A direct consequence of the above proposition is that
$q_1(u+\tfrac{1}{2},v)$, $q_1(u+\tfrac{1}{2},v)$,
$q_1(u+\tfrac{1}{2},v)$ and $q_1(u+\tfrac{1}{2},v)$ are orthogonal
to $\nu(u,v)$. We shall refer to this property by saying that {\it
crosses are planar} (see Figure \ref{figure1}). Nets with planar
crosses are called {\it asymptotic nets} (\cite{Bobenko05}). It is
also worthwhile to observe that the signs of
$(q(u+1,v+1)-q(u,v))\cdot\nu(u,v)$,
$(q(u-1,v+1)-q(u,v))\cdot\nu(u,v)$,
$(q(u-1,v-1)-q(u,v))\cdot\nu(u,v)$ and
$(q(u+1,v-1)-q(u,v))\cdot\nu(u,v)$ are alternating, and thus every
point of the surface is a {\it saddle point}.

\subsection{The affine normal map}

The affine normal map $\xi(u+\tfrac{1}{2},v+\tfrac{1}{2})$ is
defined to be
\begin{equation*}
\xi(u+\tfrac{1}{2},v+\tfrac{1}{2})=
\frac{q_{12}(u+\tfrac{1}{2},v+\tfrac{1}{2})}{F(u+\tfrac{1}{2},v+\tfrac{1}{2})}
\end{equation*}

\begin{proposition} The affine normal enjoys the following properties :

\begin{enumerate}

\item
\begin{eqnarray*}
\nu(u,v)\cdot\xi(u+\tfrac{1}{2},v+\tfrac{1}{2})&=&1\\
\nu(u+1,v)\cdot\xi(u+\tfrac{1}{2},v+\tfrac{1}{2})&=&1\\
\nu(u,v+1)\cdot\xi(u+\tfrac{1}{2},v+\tfrac{1}{2})&=&1\\
\nu(u+1,v+1)\cdot\xi(u+\tfrac{1}{2},v+\tfrac{1}{2})&=&1.
\end{eqnarray*}

\item
\begin{eqnarray*}
\nu_1(u+\tfrac{1}{2},v)\times\nu_2(u,v+\tfrac{1}{2})&=&
-F(u+\tfrac{1}{2},v+\tfrac{1}{2})\xi(u+\tfrac{1}{2},v+\tfrac{1}{2})\\
\nu_1(u+\tfrac{1}{2},v)\times\nu_2(u+1,v+\tfrac{1}{2})&=&
-F(u+\tfrac{1}{2},v+\tfrac{1}{2})\xi(u+\tfrac{1}{2},v+\tfrac{1}{2})\\
\nu_1(u+\tfrac{1}{2},v+1)\times\nu_2(u,v+\tfrac{1}{2})&=&
-F(u+\tfrac{1}{2},v+\tfrac{1}{2})\xi(u+\tfrac{1}{2},v+\tfrac{1}{2})\\
\nu_1(u+\tfrac{1}{2},v+1)\times\nu_2(u+1,v+\tfrac{1}{2})&=&
-F(u+\tfrac{1}{2},v+\tfrac{1}{2})\xi(u+\tfrac{1}{2},v+\tfrac{1}{2}).
\end{eqnarray*}



\end{enumerate}
\end{proposition}

\begin{proof}
 All formulas of Item 1 follow directly from the equation
\begin{equation*}
q_{12}(u+\tfrac{1}{2},v+\tfrac{1}{2})=\nu(u,v+1)\times\nu(u+1,v+1)-\nu(u,v)\times\nu(u+1,v).
\end{equation*}
For the second Item, we shall prove one of the equations, the
others being similar:
\begin{eqnarray*}
\nu_1(u+\tfrac{1}{2},v)\times \nu_2(u+1,v+\tfrac{1}{2})&=&
\nu(u+1,v)\times\nu(u,v+1)-\nu(u+1,v)\times\nu(u,v)-\nu(u,v)\times\nu(u,v+1)\\
&=&-q_{12}(u+\tfrac{1}{2},v+\tfrac{1}{2})\\
&=&-F(u+\tfrac{1}{2},v+\tfrac{1}{2})\xi(u+\tfrac{1}{2},v+\tfrac{1}{2}),
\end{eqnarray*}
thus proving the proposition. \end{proof}

\subsection{ Bi-linear interpolation }\label{interpolation}

The bi-linear interpolation between four points $q(u,v)$,
$q(u+1,v)$, $q(u,v+1)$ and $q(u+1,v+1)$ suits very well to the
discrete affine minimal surface with indefinite metric. This
interpolation generates a continuous surface and respects the
normal and co-normal vectors. All figures of this paper were
computed using this interpolation.

A parameterization of the hyperbolic paraboloid that passes
through $q(u,v)$, $q(u+1,v)$, $q(u,v+1)$ and $q(u+1,v+1)$ is given
by
\begin{equation}\label{patch}
r(s,t)=q(u,v)+s(q(u\!+\!1,v)-q(u,v))+t(q(u,v\!+\!1)-q(u,v))+st(q(u\!+\!1,v\!+\!1)+q(u,v)-q(u\!+\!1,v)-q(u,v\!+\!1)),
\end{equation}
for $0\leq s\leq 1$, $0\leq t\leq 1$.

\begin{lemma} The parameterization \eqref{patch} is asymptotic and the affine area of
the quadratic patch is exactly $F(u+\tfrac{1}{2},v+\tfrac{1}{2})$.
Also, $\xi(u+\tfrac{1}{2},v+\tfrac{1}{2})$ is the constant affine
normal of the surface, and  the co-normals at the corners coincide
with $\nu(u,v)$, $\nu(u+1,v)$, $\nu(u,v+1)$ and $\nu(u+1,v+1)$.
\end{lemma}
\begin{proof}
Direct calculations shows that the area element of the surface
defined by \eqref{patch} is $F(u+\tfrac{1}{2},v+\tfrac{1}{2})dsdt$
and thus its affine area is $F(u+\tfrac{1}{2},v+\tfrac{1}{2})$.
The calculation of the affine normal and the co-normals at the
corners are straightforward.
\end{proof}

\subsection{Examples}

\begin{figure}[htb]
\centering \fsep \subfigure[ Discrete helicoid. ] {
\includegraphics[width=.23
\linewidth,clip =false]{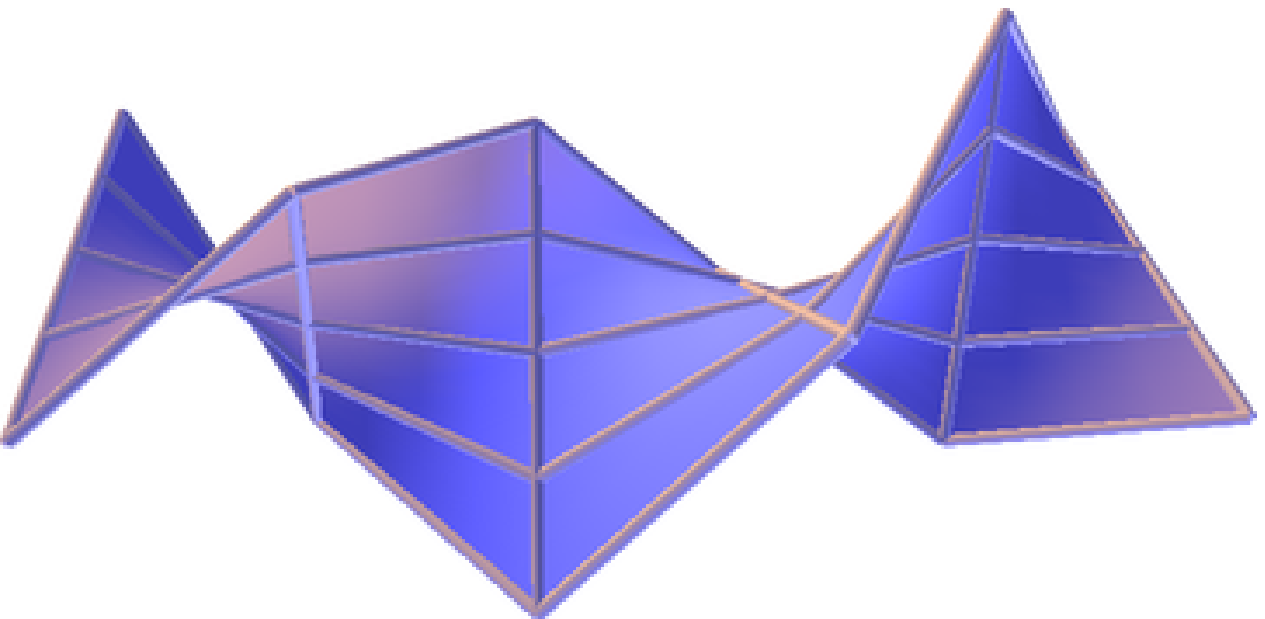}} \fsep\subfigure[
Discrete helicoid at a higher resolution. ] {
\includegraphics[width=.25\linewidth,clip
=false]{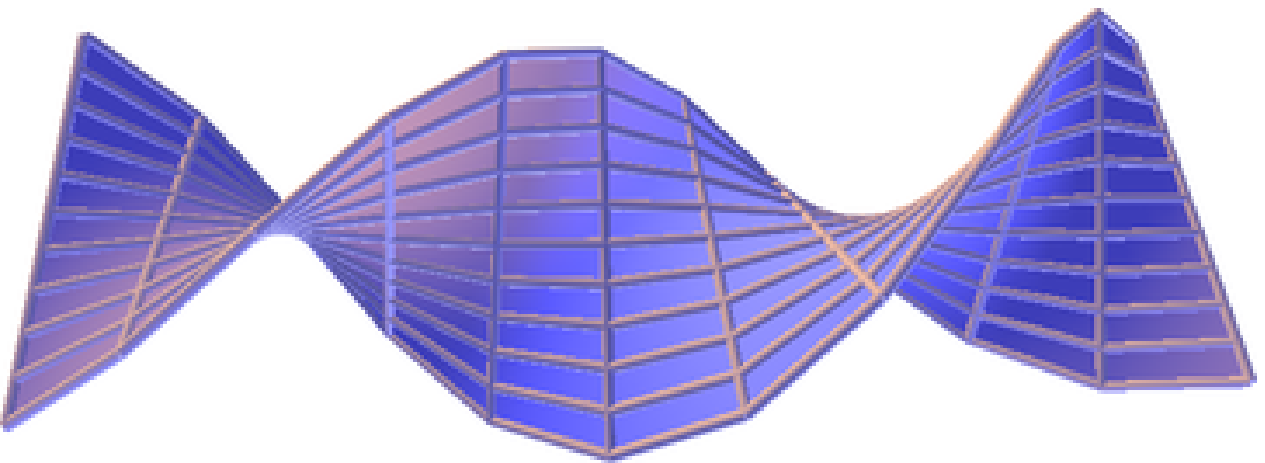}}\fsep\fsep \fsep \subfigure[ Smooth
helicoid. ] {
\includegraphics[width=.25
\linewidth,clip =false]{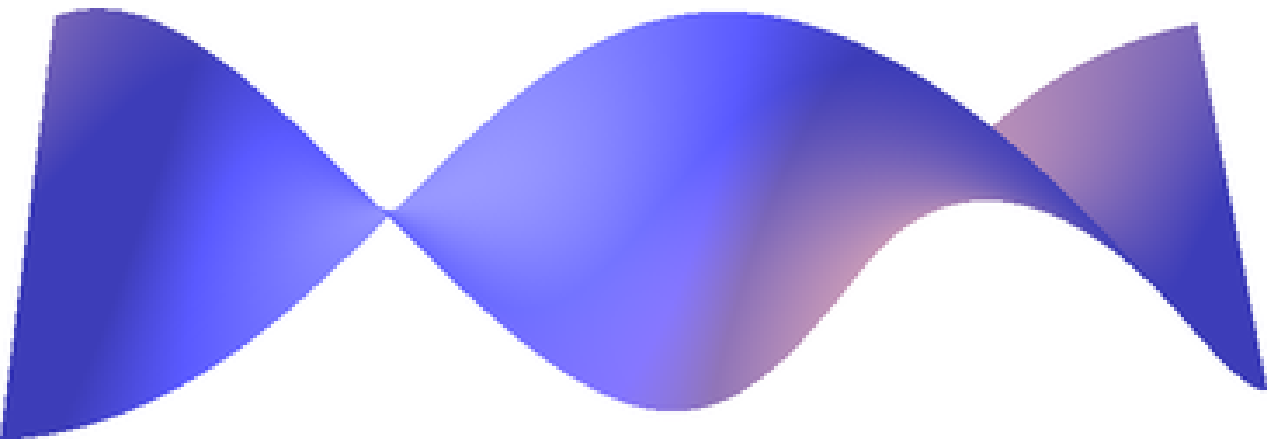}} \fsep
\caption{Discrete helicoid in two resolutions and the smooth one.}
\label{helicoid}
\end{figure}

\begin{example} The smooth helicoid can be parameterized in asymptotic
coordinates by
$$
q(u,v)=(u\cos(v),u\sin(v),v), \ (u,v)\in \R^2,
$$
and its co-normal vector field is $\nu(u,v)=(\sin(v),-\cos(v),u)$.
In order to obtain a discrete counterpart of the helicoid, we
integrate the map
$\nu(u,v)=(\sin(\frac{2\pi}{N}v),-\cos(\frac{2\pi}{N}v),u)$ for
$(u,v)\in \Z\times[0,N]\subset\Z^2$. The resulting discrete
helicoid is shown in Figure \ref{helicoid}, together with the
smooth one. We observe that the discrete parameterizations are not
restrictions to $\Z^2$ of the smooth parameterization, i.e., the
vertices of the discrete surfaces are not points of the smooth
surface.
\end{example}

\begin{figure}[htb]
\centering \fsep \subfigure[ Discrete minimal cubic. ] {
\includegraphics[width=.25
\linewidth,clip =false]{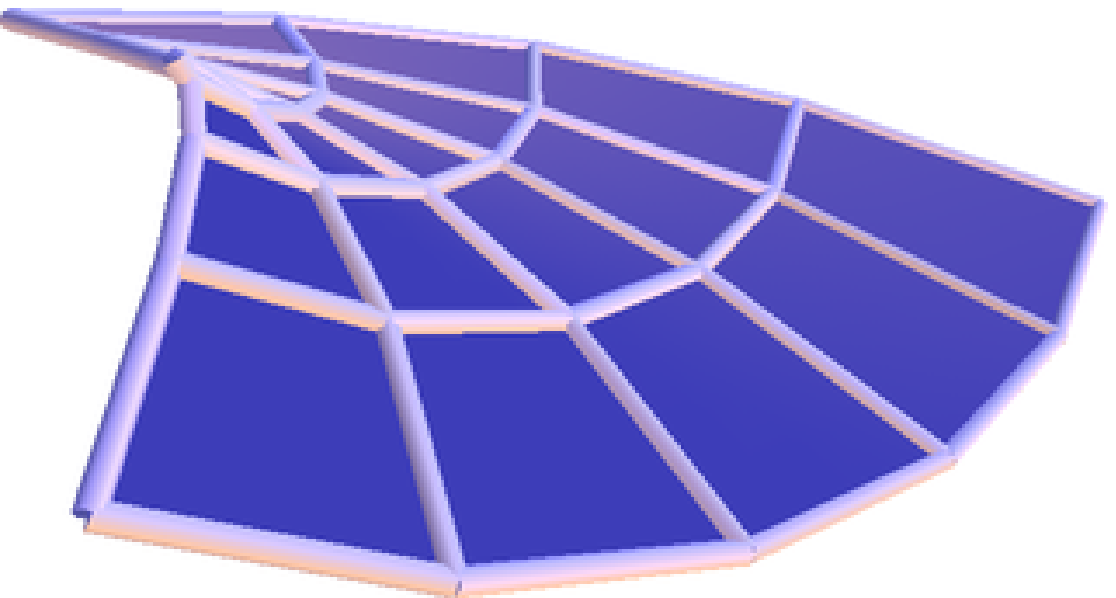}} \fsep\subfigure[
Discrete minimal cubic in higher resolution. ] {
\includegraphics[width=.25\linewidth,clip
=false]{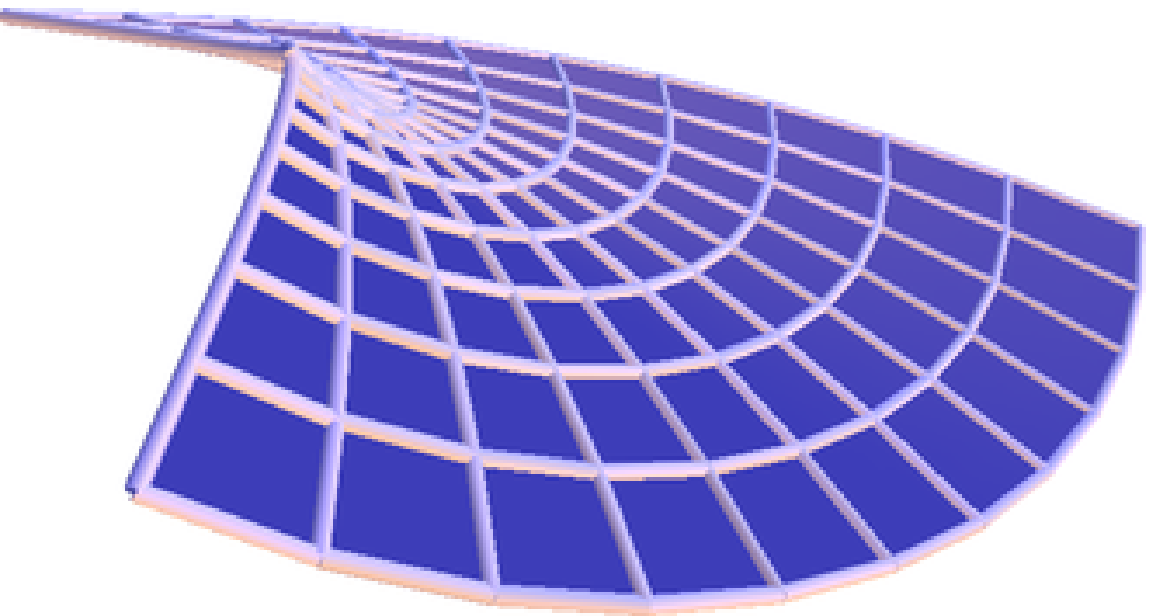}}\fsep\fsep \fsep \subfigure[ Smooth
minimal cubic. ] {
\includegraphics[width=.25
\linewidth,clip =false]{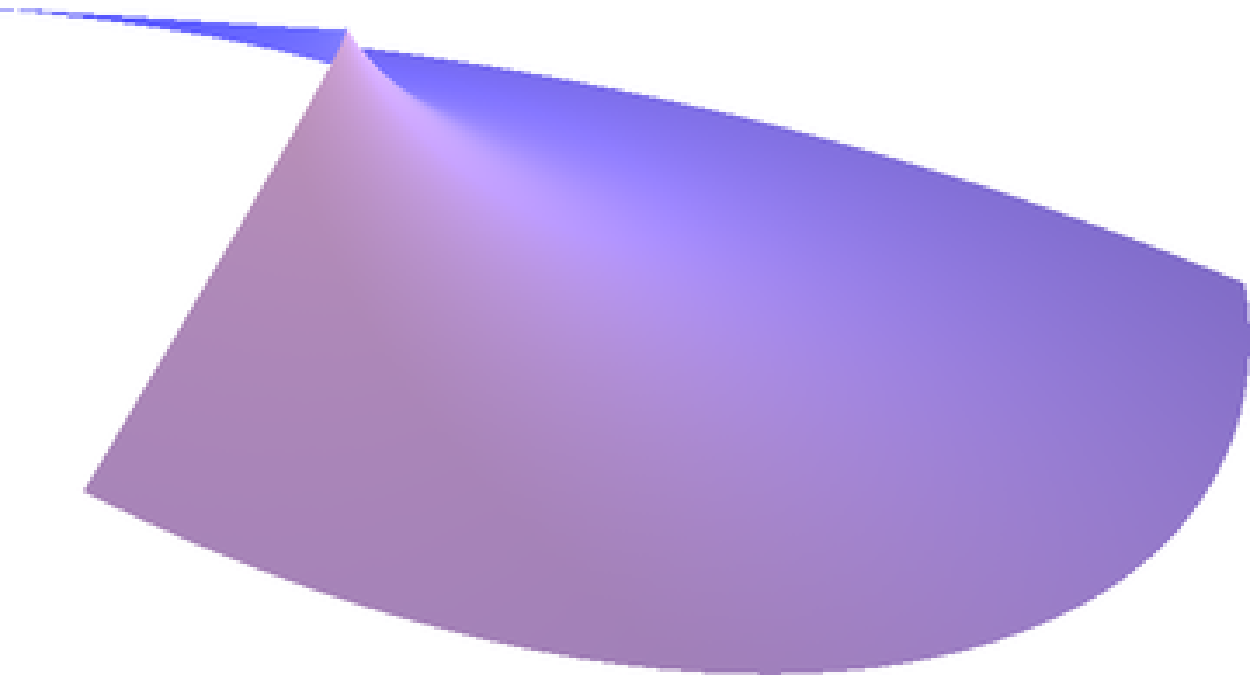}} \fsep
\caption{Discrete minimal cubic in two resolutions and smooth
minimal cubic.} \label{cubic}
\end{figure}

\begin{example} Consider a smooth vector field
$\nu(u,v)=(u,v,u^2+v^2)$, $(u,v)\in \R^2$. The associated smooth
immersion is given by
$$
q(u,v)=(u^2v-\frac{v^3}{3},v^2u-\frac{u^3}{3},-uv).
$$
To obtain the discrete counterpart of this minimal surface, we
make a discrete integration of $\nu(u,v)=(u,v,u^2+v^2)$, $(u,v)\in
\Z^2$. The resulting discrete surface, together with the smooth
one, is shown in Figure \ref{cubic}. Again, the vertices of the
discrete surface are not points of the smooth surface.
\end{example}

\begin{figure}[htb]
\centering \fsep \subfigure[ Discrete hyperbolic paraboloid. ] {
\includegraphics[width=.35
\linewidth,clip =false]{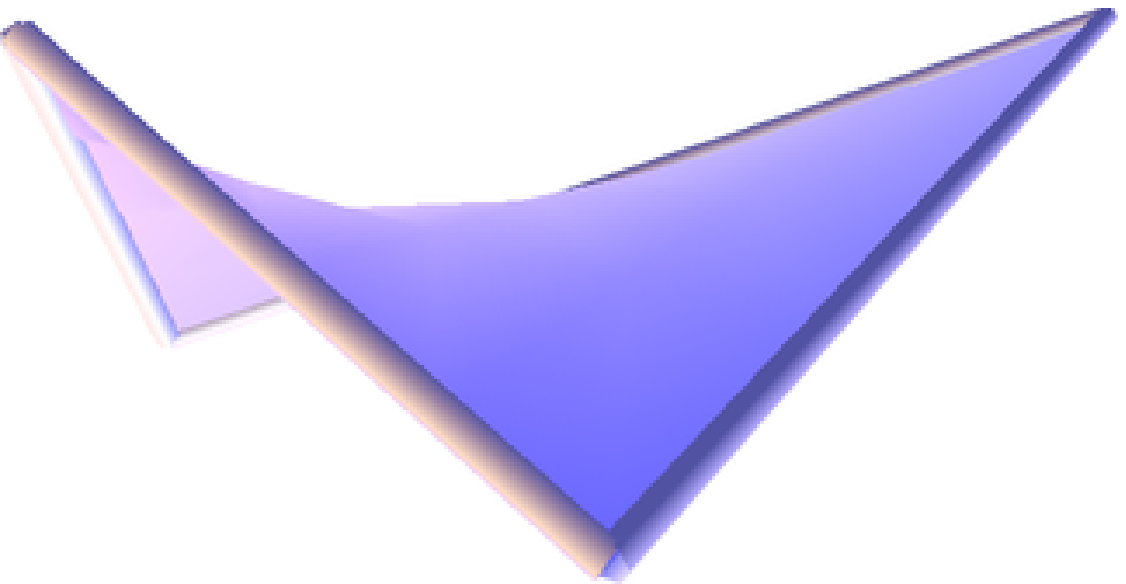}} \fsep\subfigure[
Discrete hyperbolic paraboloid in higher resolution. ] {
\includegraphics[width=.25\linewidth,clip
=false]{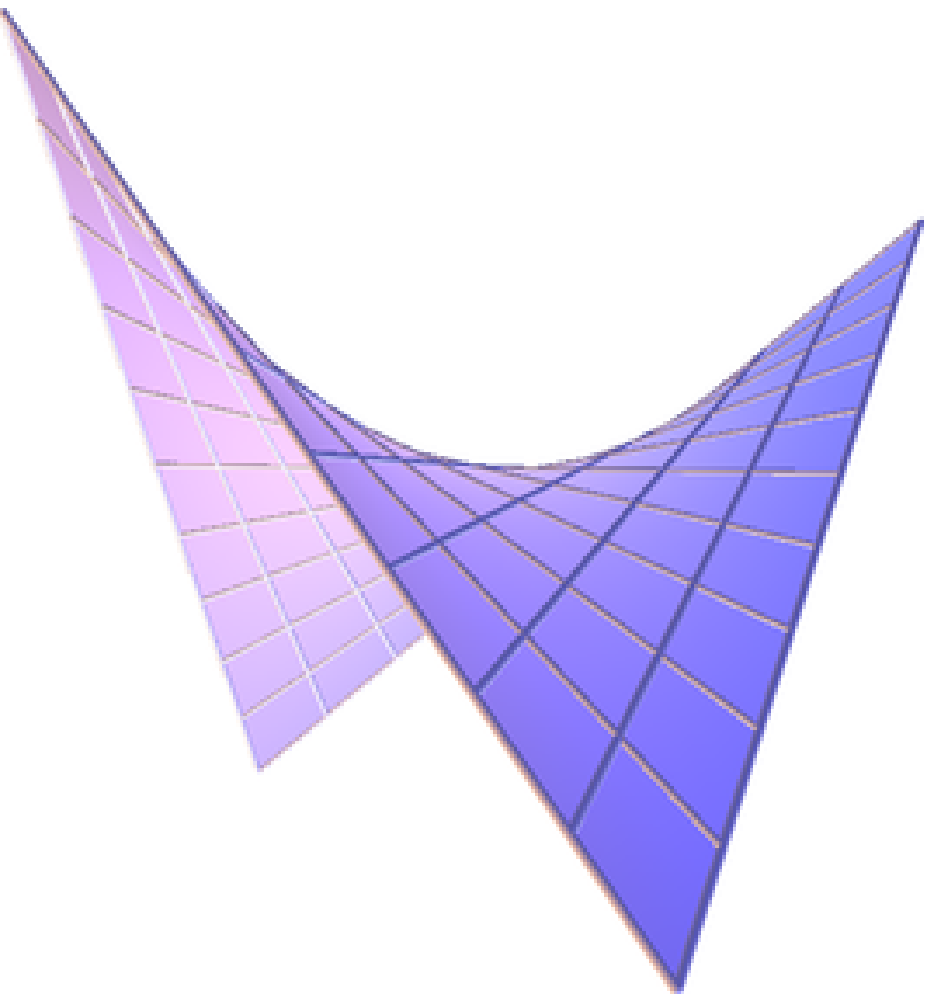}}\fsep\fsep \fsep \subfigure[ Smooth
hyperbolic paraboloid. ] {
\includegraphics[width=.25
\linewidth,clip =false]{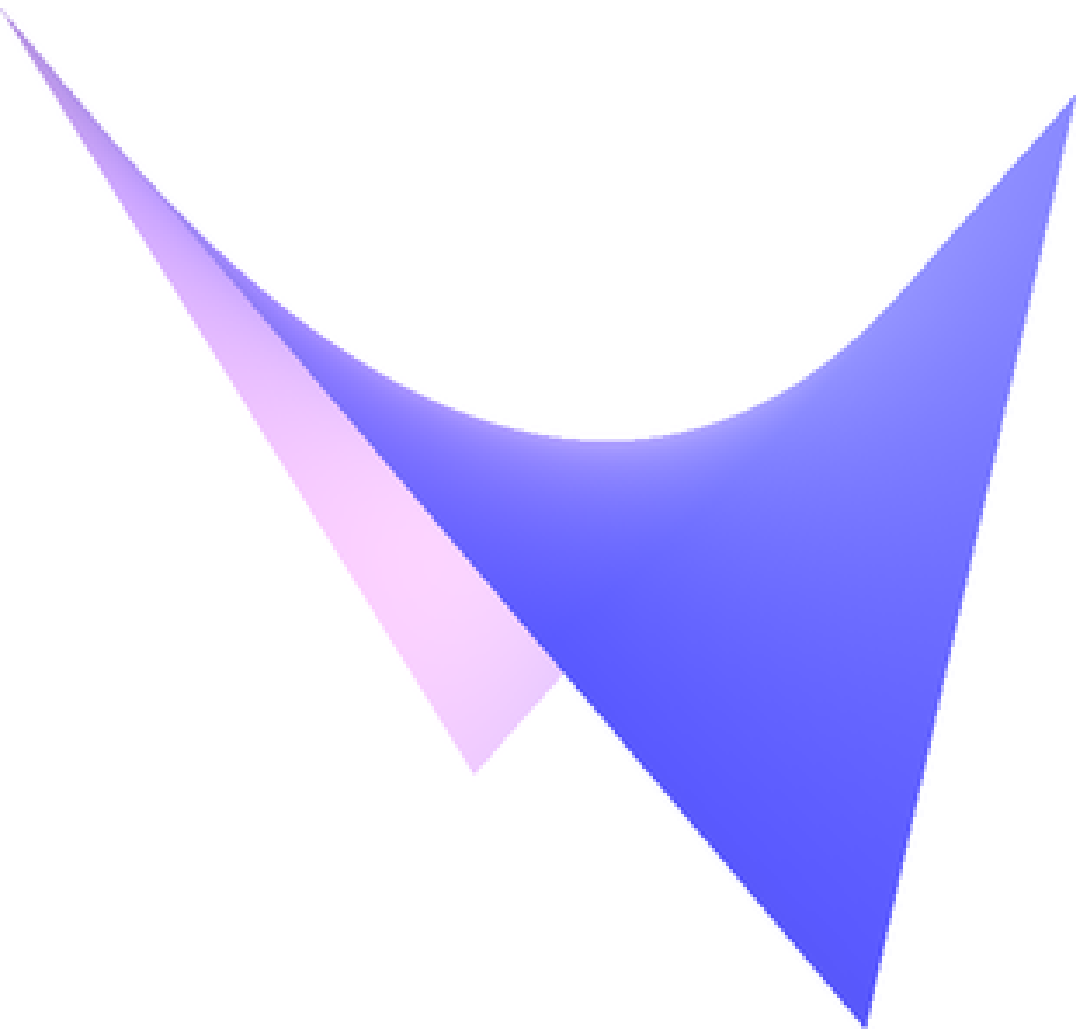}} \fsep
\caption{The discrete hyperbolic paraboloid and the smooth one
coincides.} \label{paraboloide}
\end{figure}

\begin{example} The hyperbolic paraboloid can be parameterized in asymptotic
coordinates by
$$
q(u,v)=(u,v,uv), \ (u,v)\in\R^2,
$$
and its co-normal vector field is $\nu(u,v)=(-v,-u,1)$. If we
integrate the restriction of $\nu$ to $\Z^2$, we obtain a discrete
minimal surface. It turns out that in this special case, the
discrete immersion is the restriction to $\Z^2$ of the smooth
immersion. Moreover, we observe that if we interpolate this
discrete surface as in subsection \ref{interpolation}, we obtain
again the smooth hyperbolic paraboloid (see Figure
\ref{paraboloide}).
\end{example}

A {\it discrete improper affine sphere} is a discrete minimal
surface for which the affine normal vector field is constant. It
can also be characterized by the fact that the co-normal vector
field is contained in a plane.

\begin{example}

Consider $\nu(u,v)=(\frac{v^2-u^2}{4},\frac{u-v}{2},-1)$. The
corresponding smooth affine immersion is
$$
q(u,v)=(\frac{u+v}{2},\frac{u^2+v^2}{4},\frac{(u-v)^3}{24}),
$$
and it is defined only for $u>v$. It is an improper affine sphere,
since the image of the co-normal vector field is contained in a
plane. This surface is the graph of the area distance (see
\cite{Giblin04}), a well-known concept in computer vision, to the
parabola $(t,\frac{t^2}{2})$, $t\in\R$ . The corresponding
discrete immersion is the graph of the area distance of the
polygon defined by $(t,\frac{t^2}{2})$, $t\in\Z$ (for details, see
\cite{Craizer08}). The smooth and discrete surfaces are shown in
Figure \ref{affinesphere}.

\begin{figure}[htb]
\centering \fsep \subfigure[ Discrete affine sphere. ] {
\includegraphics[width=.25
\linewidth,clip =false]{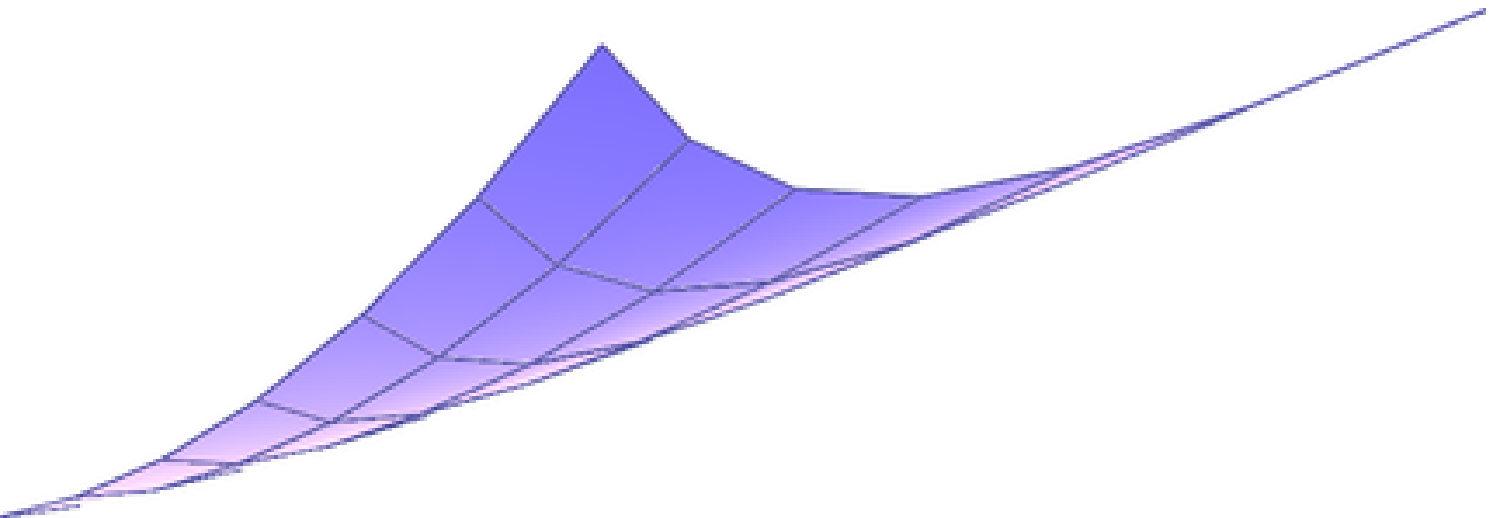}} \fsep\subfigure[
Discrete affine sphere in higher resolution. ] {
\includegraphics[width=.25\linewidth,clip
=false]{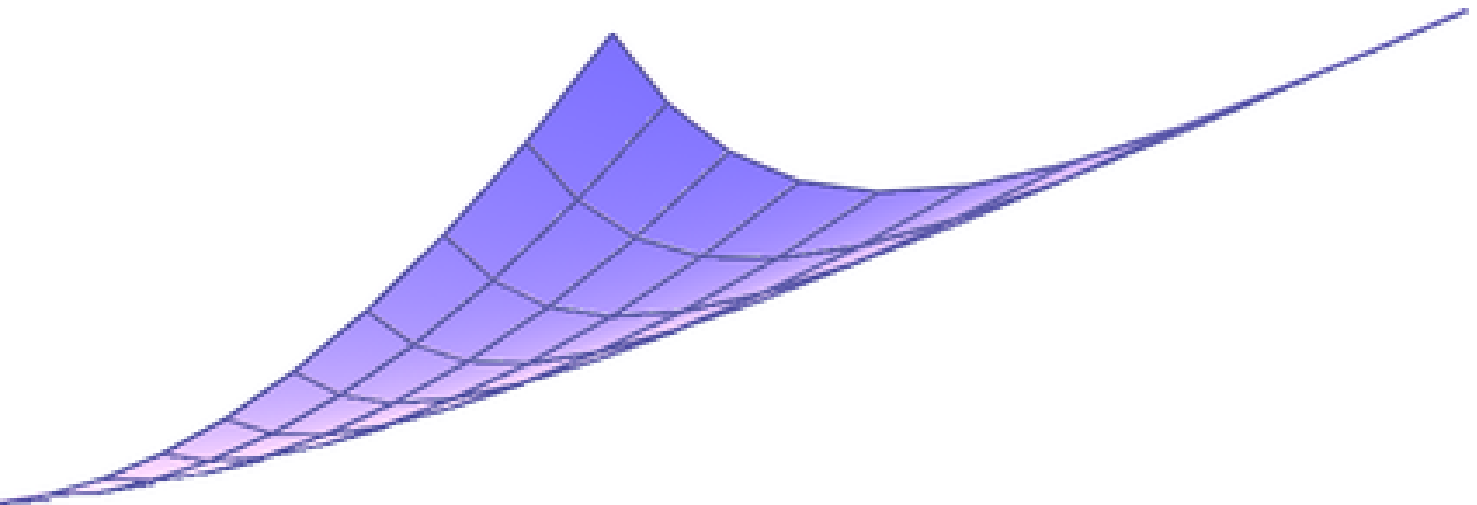}}\fsep\fsep \fsep \subfigure[ Smooth
affine sphere. ] {
\includegraphics[width=.25
\linewidth,clip =false]{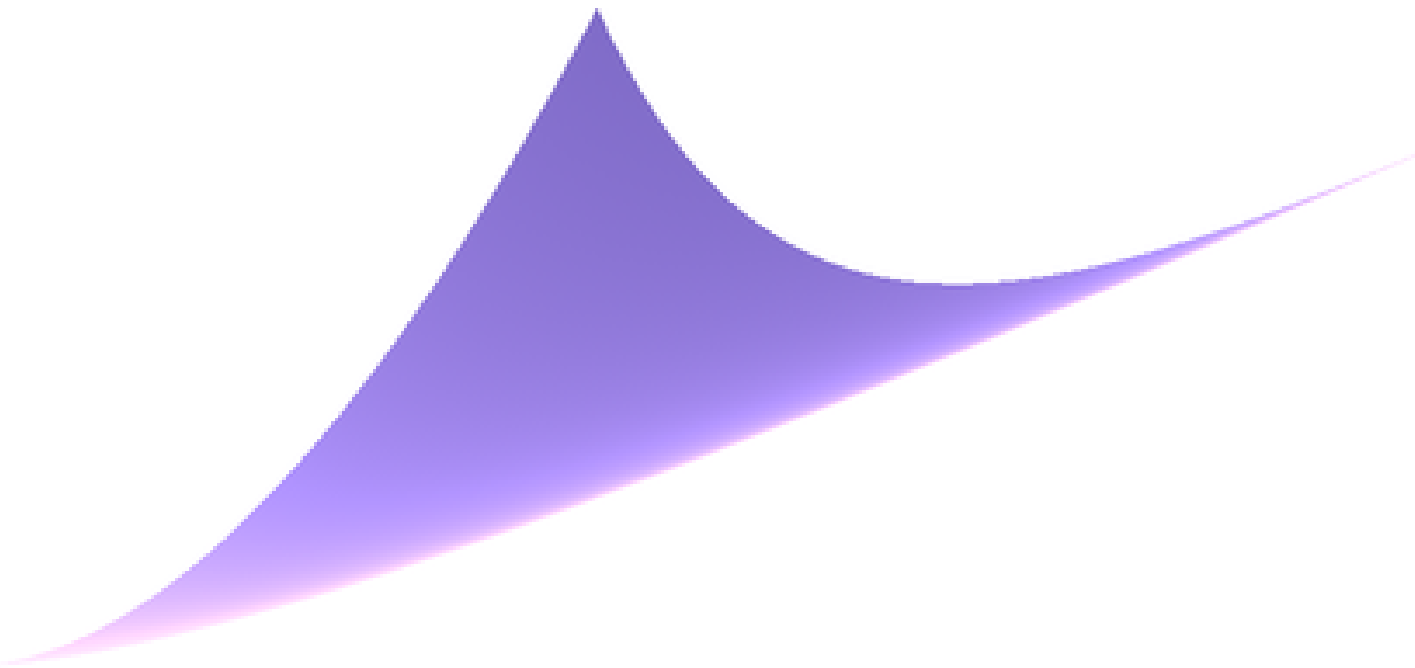}} \fsep
\caption{The discrete and smooth affine spheres.}
\label{affinesphere}
\end{figure}

\end{example}

\section{Variational property}
 In this section we introduce a functional on the space of discrete indefinite
 quadrangular surfaces
 and prove that the affine minimal discrete surfaces that we have described in Section
 \ref{SectionDefinitions} are
 actually critical points of this functional.
 \subsection{The discrete affine area functional}
 Let $S$ a discrete quadrangular surface and $q: D \to \R^3,$ with $D \subset \Z^2$ a
parameterization of $S$. We further assume that for any $ (u,v)
\in D,$ the quantity
$$ M(u,v)= \big[q(u+1,v)-q(u,v),q(u,v+1)-q(u,v),q(u+1,v+1)-q(u,v)\big]
$$ is strictly positive. The quantity $F=\sqrt{M}$ is the affine area of the
hyperbolic paraboloid that passes through the vertices
$q(u,v),q(u+1,v),q(u,v+1)$ and $q(u+1,v+1).$ The {\it affine area}
of $S$ is defined as
$$ {\cal F}(S)= \sum_{(u,v) \in D} F(u,v).$$

Let $V(u,v): D  \to \R^3$ a map such that $V(u,v)$ vanishes except
on a finite number of points $(u,v)$ of $D$.
 Intuitively,
$V$ must be regarded as a compactly supported vector field on $S$.
The surface $S(t)$ parameterized by $q_t(u,v) = q(u,v) + tV(u,v)$
is a \em deformation \em of $S.$ For $t$ small enough, we still
have $M_t(u,v) >0,$ so the next definition makes sense:

\begin{definition}
 A quadrangular surface is said to be {\it variationally discrete
affine minimal} if
$$ \left. \frac{d {\cal F}(S_t)}{dt}\right|_{t=0} =0,$$
for any such deformation.
\end{definition}


\begin{theorem}\label{variational}
Let $q: D\subset \Z^2 \to \R^3$ be a discrete affine minimal
immersion as defined in Section \ref{SectionDefinitions}. Then it
is variationally minimal.
\end{theorem}

\begin{proof} We first observe that the first
variation $\left. \frac{d {\cal F}(S_t)}{dt}\right|_{t=0}$ is
linear with respect to $V,$ so that it is sufficient to look at a
point-wise deformation. Let $(q_0,q_1,q_2,q_3)$ be a quadrangle,
whose last vertex $q_3(t)$ is deformed by
$$q_3(t)= q_3(0) +tV + o(t),$$
 i.e.  $q_3'(0)=V$. Since $F^2(t)=[q_1-q_0,q_2-q_0,q_3(t)-q_0]$,
 we obtain
$$
F'(0)=(\frac{(q_1-q_0)\times(q_2-q_0)}{2F(0)})\cdot V .
$$

If a vertex $q(u,v)$ is deformed by
$$q(u,v,t)=q(u,v)+tV+o(t),$$
it affects the affine area of its four neighbors quadrangles. The
area variation of the quadrangle $(u-\tfrac{1}{2},v-\tfrac{1}{2})$
is given by $h_1\cdot V$, where
$$
h_1 = \frac{q_1(u-\tfrac{1}{2},v-1)\times
q_2(u-1,v-\tfrac{1}{2})}{2 F(u-\tfrac{1}{2},v-\tfrac{1}{2})}
$$
Similarly, the area variations of the quadrangles
$(u+\tfrac{1}{2},v-\tfrac{1}{2})$,
$(u+\tfrac{1}{2},v+\tfrac{1}{2})$ and
$(u-\tfrac{1}{2},v+\tfrac{1}{2})$ are given by $h_2\cdot V$,
$h_3\cdot V$ and $h_4\cdot V$, where
\begin{eqnarray*}
h_2&=&-\frac{q_1(u+\tfrac{1}{2},v-1)\times
q_2(u+1,v-\tfrac{1}{2})}{2F(u+\tfrac{1}{2},v-\tfrac{1}{2})}\\
h_3&=&\frac{q_1(u+\tfrac{1}{2},v+1)\times
q_2(u+1,v+\tfrac{1}{2})}{2F(u+\tfrac{1}{2},v+\tfrac{1}{2})}\\
h_4&=&-\frac{q_1(u-\tfrac{1}{2},v+1)\times
q_2(u-1,v+\tfrac{1}{2})}{2F(u-\tfrac{1}{2},v+\tfrac{1}{2})}.
\end{eqnarray*}
Since $\left. \frac{d {\cal
F}(S_t)}{dt}\right|_{t=0}=(h_1+h_2+h_3+h_4)\cdot V$, the surface
is variationally minimal if and only if $h_1+h_2+h_3+h_4=0$, for
any $(u,v)\in D$.

Assuming that $S$ is affine minimal, we have that
$$
\nu(u+1,v+1)+\nu(u,v)-\nu(u,v+1)-\nu(u+1,v)=0,
$$
for any $(u,v)\in D$, implying that
$$
\nu(u-1,v-1)+\nu(u+1,v+1)-\nu(u-1,v+1)-\nu(u+1,v-1)=0,
$$
for any $(u,v)\in D$, which, by Proposition \ref{DefMinimas}, is
equivalent to $h_1+h_2+h_3+h_4=0$.
\end{proof}

\section{Structural equations and compatibility}

In this section we define the {\it discrete affine cubic form} and
show that any discrete affine minimal surface must satisfy
compatibility equations that involve also the {\it discrete
quadratic form}, i.e., the Berwald-Blaschke metric. On the other
hand, given discrete quadratic and cubic forms satisfying the
compatibility equations, there exists a discrete affine minimal
surface, unique up to affine transformations of $\R^3$, with the
given quadratic and cubic forms. This result is the discrete
counterpart of the structural theorem for smooth affine minimal
surfaces.

\subsection{The discrete cubic form}\label{SectionCubic}

We define the discrete cubic form as $A(u,v)\delta
u^3+B(u,v)\delta v^3$, where
\begin{eqnarray*}
A(u,v)&=&[q_{1}(u-\tfrac{1}{2},v),q_1(u+\tfrac{1}{2},v),\xi(u\pm\tfrac{1}{2},v\pm\tfrac{1}{2})]\\
B(u,v)&=&[q_{2}(u,v+\tfrac{1}{2}),q_2(u,v-\tfrac{1}{2}),\xi(u\pm\tfrac{1}{2},v\pm\tfrac{1}{2})].
\end{eqnarray*}
Since we are interested only in the coefficients $A(u,v)$ and
$B(u,v)$ of the discrete cubic form, we shall not discuss in this
paper the meaning of the symbols $\delta u^3$ and $\delta v^3$.

From the definition of $A$ and $B$, we can write
\begin{eqnarray*}
q_{11}(u,v)&=&\frac{1}{F(u+\tfrac{1}{2},v+\tfrac{1}{2})}
(F_1(u,v+\tfrac{1}{2})q_1(u+\tfrac{1}{2},v)+A(u,v)q_2(u,v+\tfrac{1}{2}))\\
&=&\frac{1}{F(u-\tfrac{1}{2},v+\tfrac{1}{2})}
(F_1(u,v+\tfrac{1}{2})q_1(u-\tfrac{1}{2},v)+A(u,v)q_2(u,v+\tfrac{1}{2}))\\
&=&\frac{1}{F(u+\tfrac{1}{2},v-\tfrac{1}{2})}
(F_1(u,v-\tfrac{1}{2})q_1(u+\tfrac{1}{2},v)+A(u,v)q_2(u,v-\tfrac{1}{2}))\\
&=&\frac{1}{F(u-\tfrac{1}{2},v-\tfrac{1}{2})}
(F_1(u,v-\tfrac{1}{2})q_1(u-\tfrac{1}{2},v)+A(u,v)q_2(u,v-\tfrac{1}{2}))\\
q_{22}(u,v)&=&\frac{1}{F(u+\tfrac{1}{2},v+\tfrac{1}{2})}
(B(u,v)q_1(u+\tfrac{1}{2},v)+F_2(u+\tfrac{1}{2},v)q_2(u,v+\tfrac{1}{2}))\\
&=&\frac{1}{F(u-\tfrac{1}{2},v+\tfrac{1}{2})}
(B(u,v)q_1(u-\tfrac{1}{2},v)+F_2(u-\tfrac{1}{2},v)q_2(u,v+\tfrac{1}{2}))\\
&=&\frac{1}{F(u+\tfrac{1}{2},v-\tfrac{1}{2})}(
B(u,v)q_1(u+\tfrac{1}{2},v)+F_2(u+\tfrac{1}{2},v)q_2(u,v-\tfrac{1}{2}))\\
&=&\frac{1}{F(u-\tfrac{1}{2},v-\tfrac{1}{2})}(
B(u,v)q_1(u-\tfrac{1}{2},v)+F_2(u-\tfrac{1}{2},v)q_2(u,v-\tfrac{1}{2})),
\end{eqnarray*}
where
$F_{1}(u,v+\tfrac{1}{2})=F(u+\tfrac{1}{2},v+\tfrac{1}{2})-F(u-\tfrac{1}{2},v+\tfrac{1}{2})$
and
$F_{2}(u+\tfrac{1}{2},v)=F(u+\tfrac{1}{2},v+\tfrac{1}{2})-F(u+\tfrac{1}{2},v-\tfrac{1}{2})$.

\subsection{ Derivatives of the affine normal }

We shall now calculate the derivatives of the affine normal. We
first prove a technical lemma:

\begin{lemma}\label{technical}
The discrete derivatives $A_2$ and $B_1$ can be expressed as:
\begin{eqnarray*}
A_2(u,v+\tfrac{1}{2})&=&F(u-\tfrac{1}{2},v+\tfrac{1}{2})
[q_1(u+\tfrac{1}{2},v),\xi(u-\tfrac{1}{2},v+\tfrac{1}{2}),\xi(u+\tfrac{1}{2},v+\tfrac{1}{2})]\\
B_1(u+\tfrac{1}{2},v)&=&-F(u+\tfrac{1}{2},v-\tfrac{1}{2})
[q_2(u,v+\tfrac{1}{2}),\xi(u+\tfrac{1}{2},v-\tfrac{1}{2}),\xi(u+\tfrac{1}{2},v+\tfrac{1}{2})].
\end{eqnarray*}
\end{lemma}
\begin{proof} We can write
$$
q_1(u-\tfrac{1}{2},v)\times q_1(u+\tfrac{1}{2},v)=A(u,v)\nu(u,v)
$$
Differentiating with respect to $v$ we obtain
$$
A(u,v+1)\nu(u,v+1)-A(u,v)\nu(u,v) = q_1(u-\tfrac{1}{2},v+1)\times
q_1(u+\tfrac{1}{2},v+1)-q_1(u-\tfrac{1}{2},v)\times
q_1(u+\tfrac{1}{2},v)
$$
Multiplying by $\xi(u+\tfrac{1}{2},v+\tfrac{1}{2})$ we have
\begin{eqnarray*}
A(u,v+1)-A(u,v)&=&[q_1(u-\tfrac{1}{2},v+1)\times
q_1(u+\tfrac{1}{2},v+1)-q_1(u-\tfrac{1}{2},v)\times
q_1(u+\tfrac{1}{2},v)]\cdot \xi(u+\tfrac{1}{2},v+\tfrac{1}{2})\\
&=&
[q_1(u-\tfrac{1}{2},v+1)-q_1(u-\tfrac{1}{2},v),q_1(u+\tfrac{1}{2},v),
\xi(u+\tfrac{1}{2},v+\tfrac{1}{2})]\\
&=& -F(u-\tfrac{1}{2},v+\tfrac{1}{2})
[\xi(u-\tfrac{1}{2},v+\tfrac{1}{2}),q_1(u+\tfrac{1}{2},v),\xi(u+\tfrac{1}{2},v+\tfrac{1}{2})]
\end{eqnarray*}
The calculation for $B_1$ is similar. \end{proof}
\bigskip

Observe that
$\nu(u,v)\cdot\xi_1(u,v+\tfrac{1}{2})=\nu(u,v)\cdot\xi_2(u+\tfrac{1}{2},v)=0$,
and thus we can write $\xi_1(u,v+\tfrac{1}{2})$ and
$\xi_2(u+\tfrac{1}{2},v)$ as linear combinations of
$q_1(u+\tfrac{1}{2},v)$ and $q_2(u,v+\tfrac{1}{2}))$. More
precisely, we have the following proposition:

\begin{proposition}
The discrete derivatives of the affine normals can be expressed
as:

\begin{eqnarray}
F(u-\tfrac{1}{2},v+\tfrac{1}{2})F(u+\tfrac{1}{2},v+\tfrac{1}{2})\xi_{1}(u,v+\tfrac{1}{2})&=&
-A_2(u,v+\tfrac{1}{2})q_2(u,v+\tfrac{1}{2})\label{DerivativesNormala}\\
F(u+\tfrac{1}{2},v-\tfrac{1}{2})F(u+\tfrac{1}{2},v+\tfrac{1}{2})\xi_{2}(u+\tfrac{1}{2},v)&=&
-B_1(u+\tfrac{1}{2},v)q_1(u+\tfrac{1}{2},v).\label{DerivativesNormalb}
\end{eqnarray}

\end{proposition}

\begin{proof} We first show that the coefficient of
$q_1(u+\tfrac{1}{2},v)$ in the expansion of
$\xi_{1}(u,v+\tfrac{1}{2})$ is zero. We have
\begin{eqnarray*}
[\xi_{1}(u,v+\tfrac{1}{2}),q_2(u,v+\tfrac{1}{2}),\xi(u+\tfrac{1}{2},v+\tfrac{1}{2})]&=&
-[\xi(u-\tfrac{1}{2},v+\tfrac{1}{2}),q_2(u,v+\tfrac{1}{2}),\xi(u+\tfrac{1}{2},v+\tfrac{1}{2})]\\
&=&
\frac{[q_2(u-1,v+\tfrac{1}{2}),q_2(u,v+\tfrac{1}{2}),q_2(u+1,v+\tfrac{1}{2})]}
{F(u-\tfrac{1}{2},v+\tfrac{1}{2})F(u+\tfrac{1}{2},v+\tfrac{1}{2})}
\end{eqnarray*}
And
\begin{eqnarray*}
q_2(u-1,v+\tfrac{1}{2})\times q_2(u,v+\tfrac{1}{2})&=&
(\nu(u-1,v+1)\times \nu(u-1,v))\times (\nu(u,v+1)\times
\nu(u,v))\\
&=& ((\nu(u,v+1)-\nu(u,v))\times \nu(u-1,v))\times
(\nu(u,v+1)\times \nu(u,v))\\
&=& -F(u-\tfrac{1}{2},v+\tfrac{1}{2})(\nu(u,v+1)-\nu(u,v))
\end{eqnarray*}
So
\begin{eqnarray*}
[\xi_{1}(u,v+\tfrac{1}{2}),q_2(u,v+\tfrac{1}{2}),\xi(u+\tfrac{1}{2},v+\tfrac{1}{2})]&=&
-\frac{(\nu(u,v+1)-\nu(u,v))\cdot (\nu(u+1,v+1)\times
\nu(u+1,v))}{F(u+\tfrac{1}{2},v+\tfrac{1}{2})}
\\
&=&
-\frac{1}{F(u+\tfrac{1}{2},v+\tfrac{1}{2})}(F(u+\tfrac{1}{2},v+\tfrac{1}{2})-F(u+\tfrac{1}{2},v+\tfrac{1}{2}))=0.
\end{eqnarray*}
We can now easily complete the proof of the first equation using
lemma \ref{technical}. The proof of the second equation is
similar.
\end{proof}

\begin{corollary} A discrete affine minimal surface is an improper affine sphere if and only if $A=A(u)$ and
$B=B(v)$.
\end{corollary}

\subsection{Compatibility equations}

In this subsection we obtain three compatibility equations. They
are generalizations of the equations obtained in \cite{Matsuura03}
for discrete improper affine spheres. The first equation is proved
in the following lemma:

\begin{lemma}
\begin{equation}\label{Compatibility0}
F(u-\tfrac{1}{2},v+\tfrac{1}{2})F(u+\tfrac{1}{2},v-\tfrac{1}{2})-
F(u+\tfrac{1}{2},v+\tfrac{1}{2})F(u-\tfrac{1}{2},v-\tfrac{1}{2})=A(u,v)B(u,v).
\end{equation}
\end{lemma}
\begin{proof}
We can calculate $q_{112}(u,v+\tfrac{1}{2})$ as
$q_{12}(u+\tfrac{1}{2},v+\tfrac{1}{2})-q_{12}(u-\tfrac{1}{2},v+\tfrac{1}{2})$
and also as $q_{11}(u,v+1)-q_{11}(u,v)$. Calculating in the first
way, we have from $q_{12}=F\xi$ that
\begin{eqnarray*}
q_{112}(u,v+\tfrac{1}{2})&=&
F_1(u,v+\tfrac{1}{2})\xi(u+\tfrac{1}{2},v+\tfrac{1}{2})
+F(u+\tfrac{1}{2},v+\tfrac{1}{2})\xi_1(u,v+\tfrac{1}{2})\\
&=& F_1(u,v+\tfrac{1}{2})\xi(u+\tfrac{1}{2},v+\tfrac{1}{2})
-\frac{A_2(u,v+\tfrac{1}{2})}{F(u-\tfrac{1}{2},v+\tfrac{1}{2})}q_2(u,v+\tfrac{1}{2}).
\end{eqnarray*}
Calculating in the second way, formulas of subsection
\ref{SectionCubic} imply that
\begin{eqnarray*}
q_{112}(u,v+\tfrac{1}{2})&=&
\left(\frac{F_1(u,v+\tfrac{1}{2})}{F(u+\tfrac{1}{2},v+\tfrac{1}{2})}-\frac{F_1(u,v-\tfrac{1}{2})}{F(u+\tfrac{1}{2},v-\tfrac{1}{2})}\right)
q_1(u+\tfrac{1}{2},v)
+F_1(u,v+\tfrac{1}{2})\xi(u+\tfrac{1}{2},v+\tfrac{1}{2})\\
&+&
\left(\frac{A(u,v+1)}{F(u+\tfrac{1}{2},v+\tfrac{1}{2})}-\frac{A(u,v)}{F(u+\tfrac{1}{2},v-\tfrac{1}{2})}\right)q_2(u,v+\tfrac{1}{2})
+\frac{A(u,v)}{F(u+\tfrac{1}{2},v-\tfrac{1}{2})}q_{22}(u,v)
\end{eqnarray*}
Now, using the formula for $q_{22}(u,v)$ and comparing the
coefficients of $q_1(u+\tfrac{1}{2},v)$, we obtain
\begin{equation*}
A(u,v)B(u,v)+F_1(u,v+\tfrac{1}{2})F(u+\tfrac{1}{2},v-\tfrac{1}{2})-F_1(u,v-\tfrac{1}{2})F(u+\tfrac{1}{2},v+\tfrac{1}{2})=0,
\end{equation*}
thus proving the lemma.

\end{proof}

\medskip
The two other compatibility equations are obtained from Equations
\eqref{DerivativesNormala} and \eqref{DerivativesNormalb}. We can
write
\begin{eqnarray*}
F(u+\tfrac{1}{2},v+\tfrac{1}{2})\xi_{12}(u,v)&=&
-\frac{B(u,v)A_2(u,v-\tfrac{1}{2})}{F(u-\tfrac{1}{2},v-\tfrac{1}{2})
F(u+\tfrac{1}{2},v-\tfrac{1}{2})}q_1(u+\tfrac{1}{2},v)\\
&+&\left(\frac{A_2(u,v-\tfrac{1}{2})}{F(u-\tfrac{1}{2},v-\tfrac{1}{2})}-\frac{A_2(u,v+\tfrac{1}{2})}
{F(u-\tfrac{1}{2},v+\tfrac{1}{2})}\right)q_2(u,v+\tfrac{1}{2})\\
F(u+\tfrac{1}{2},v+\tfrac{1}{2})\xi_{21}(u,v)&=&
\left(\frac{B_1(u-\tfrac{1}{2},v)}{F(u-\tfrac{1}{2},v-\tfrac{1}{2})}-\frac{B_1(u+\tfrac{1}{2},v)}
{F(u+\tfrac{1}{2},v-\tfrac{1}{2})}\right)q_1(u+\tfrac{1}{2},v)\\
&-&\frac{A(u,v)B_1(u-\tfrac{1}{2},v)}{F(u-\tfrac{1}{2},v-\tfrac{1}{2})
F(u-\tfrac{1}{2},v+\tfrac{1}{2})}q_2(u,v+\tfrac{1}{2}).
\end{eqnarray*}
Thus we get
\begin{eqnarray}
F(u-\tfrac{1}{2},v-\tfrac{1}{2})B_1(u+\tfrac{1}{2},v)-F(u+\tfrac{1}{2},v-\tfrac{1}{2})B_1(u-\tfrac{1}{2},v)&=&
B(u,v)A_2(u,v-\tfrac{1}{2})\label{CompatibilityDisca}\\
F(u-\tfrac{1}{2},v-\tfrac{1}{2})A_2(u,v+\tfrac{1}{2})-F(u-\tfrac{1}{2},v+\tfrac{1}{2})A_2(u,v-\tfrac{1}{2})&=&
A(u,v)B_1(u-\tfrac{1}{2},v).\label{CompatibilityDiscb}
\end{eqnarray}

\subsection{Existence and uniqueness theorem}

In this subsection, we prove the existence and uniqueness of a
discrete affine minimal surface with given quadratic and cubic
forms satisfying the compatibility equations.

\begin{theorem}\label{Existence}
Given function $F(u+\tfrac{1}{2},v+\tfrac{1}{2})$, $A(u,v)$ and
$B(u,v)$ satisfying the compatibility equations
\eqref{Compatibility0}, \eqref{CompatibilityDisca} and
\eqref{CompatibilityDiscb}, there exists a discrete affine minimal
surface $q(u,v)$ with quadratic form $Fdudv$ and cubic form
$A\delta u^3+B\delta v^3$. Moreover, two discrete affine minimal
surfaces with the same quadratic and cubic forms are affine
equivalent.
\end{theorem}

\begin{proof}
We begin by choosing four points $q(0,0)$, $q(1,0)$, $q(0,1)$ and
$q(1,1)$ satisfying
$[q(1,0)-q(0,0),q(0,1)-q(0,0),q(1,1)-q(0,0)]=F^2(\tfrac{1}{2},\tfrac{1}{2})$.
This four points are determined up to an affine transformation of
$\R^3$.

From a quadrangle $(u-\tfrac{1}{2},v-\tfrac{1}{2})$, one can
extend the definition of $q$ to the quadrangles
$(u+\tfrac{1}{2},v-\tfrac{1}{2})$ and
$(u-\tfrac{1}{2},v+\tfrac{1}{2})$ by the formulas of Section
\ref{SectionCubic}. With these extensions, we can calculate
$\xi(u+\tfrac{1}{2},v-\tfrac{1}{2})$ and
$\xi(u-\tfrac{1}{2},v+\tfrac{1}{2})$. It is clear that
$\xi_1(u,v-\tfrac{1}{2})$ and $\xi_2(u-\tfrac{1}{2},v)$ satisfy
equations \eqref{DerivativesNormala} and
\eqref{DerivativesNormalb}. The coherence of these extensions are
assured by formula \eqref{Compatibility0}.


 Then one can extend the definition
of $q$ to $(u+\tfrac{1}{2},v+\tfrac{1}{2})$ in two different ways:
from the quadrangle $(u+\tfrac{1}{2},v-\tfrac{1}{2})$ and from the
$(u-\tfrac{1}{2},v+\tfrac{1}{2})$. Our task is to show that both
extensions leads to the same result. This amounts to check that
both affine normals $\xi(u+\tfrac{1}{2},v+\tfrac{1}{2})$ are the
same, which in fact reduces to verify that $\xi_{12}=\xi_{21}$.
But this last equation holds by the compatibility hypothesis,
which completes the proof of the theorem.
\end{proof}


\smallskip\noindent{\bf Acknowledgements.} The first and third
authors want to thank CNPq for financial support during the
preparation of this paper.


\bibliography{AffineMinimalv2}

\begin{thebibliography}{10}

\bibitem{Bobenko08}
A.~Bobenko, P.~Schr{\"o}der, J.~Sullivan, and G.~Ziegler, editors.
\newblock {\em Discrete Differential Geometry}, volume~38 of {\em Oberwolfach
  Seminars}.
\newblock Birkhauser, 2008.

\bibitem{Bobenko05}
A.~I. Bobenko and Y.~B.Suris.
\newblock Discrete differential geometry: Consistency as integrability.
\newblock {\em pre-print}, 2005.

\bibitem{Bobenko06}
A.~I. Bobenko, T.~Hoffmann, and B.~A. Springborn.
\newblock Minimal surfaces from circle patterns: Geometry from combinatorics.
\newblock {\em Annals of Mathematics}, 164(1):231--264, 2006.

\bibitem{Bobenko199}
A.~I. Bobenko and W.~K. Schief.
\newblock Affine spheres: Discretization via duality relations.
\newblock {\em Experimental Mathematics}, 8(3):261--280, 1999.

\bibitem{Buchin83}
S.~Buchin.
\newblock {\em Affine Differential Geometry}.
\newblock Science Press, Beijing, China, Gordon and Breach,Science Publishers,
  New York, 1983.

\bibitem{Craizer08}
M.~Craizer, M.~A. da~Silva, and R.~C. Teixeira.
\newblock Area distances of convex plane curves and improper affine spheres.
\newblock {\em pre-print}, 2008.

\bibitem{Calabi82}
E.Calabi.
\newblock Hypersurfaces with maximal affinely invariant area.
\newblock {\em American Journal of Mathematics}, 104:91--126, 1982.

\bibitem{Calabi88}
E.Calabi.
\newblock Convex affine maximal surfaces.
\newblock {\em Results in Mathematics}, 13:199--223, 1988.

\bibitem{Simon93}
A.-M. Li, U.~Simon, and G.~Zhao.
\newblock {\em Global Affine Differential Geometry of Hypersurfaces}.
\newblock De Gruyter Expositions in Mathematics, 1993.

\bibitem{Matsuura03}
N.~Matsuura and H.~Urakawa.
\newblock Discrete improper affine spheres.
\newblock {\em Journal of Geometry and Physics}, 45:164--183, 2003.

\bibitem{Giblin04}
M.~Niethammer, S.~Betelu, G.~Sapiro, A.~Tannenbaum, and P.~J. Giblin.
\newblock Area-based medial axis of planar curves.
\newblock {\em International Journal of Computer Vision}, 60(3):203--224, 2004.

\bibitem{Nomizu94}
K.~Nomizu and T.~Sasaki.
\newblock {\em Affine Differential Geometry}.
\newblock Cambridge University Press, 1994.

\bibitem{Vrancken89}
L.~Verstraelen and L.~Vrancken.
\newblock Affine variation formulas and affine minimal surfaces.
\newblock {\em Michigan Mathematical Journal}, 36:77--93, 1989.

\end{thebibliography}


\end{document}